\documentclass[11pt]{amsart}
\usepackage{graphicx}
\usepackage{amssymb, amsmath, amsthm, mathdots, hyperref}
\usepackage{youngtab}
\usepackage[all]{xy}
\begin{document}
\newtheorem{defn0}{Definition}[section]
\newtheorem{prop0}[defn0]{Proposition}
\newtheorem{thm0}[defn0]{Theorem}
\newtheorem{lemma0}[defn0]{Lemma}
\newtheorem{coro0}[defn0]{Corollary}
\newtheorem{exa}[defn0]{Example}
\newtheorem{rem0}[defn0]{Remark}
\newtheorem{quest0}[defn0]{Question}
\numberwithin{equation}{section}
\def\rig#1{\smash{ \mathop{\longrightarrow}
    \limits^{#1}}}
\def\nwar#1{\nwarrow
   \rlap{$\vcenter{\hbox{$\scriptstyle#1$}}$}}
\def\near#1{\nearrow
   \rlap{$\vcenter{\hbox{$\scriptstyle#1$}}$}}
\def\sear#1{\searrow
   \rlap{$\vcenter{\hbox{$\scriptstyle#1$}}$}}
\def\swar#1{\swarrow
   \rlap{$\vcenter{\hbox{$\scriptstyle#1$}}$}}
\def\dow#1{\Big\downarrow
   \rlap{$\vcenter{\hbox{$\scriptstyle#1$}}$}}
\def\up#1{\Big\uparrow
   \rlap{$\vcenter{\hbox{$\scriptstyle#1$}}$}}
\def\lef#1{\smash{ \mathop{\longleftarrow}
    \limits^{#1}}}
 \def\stackbelow#1#2{\underset{\displaystyle\overset{\displaystyle\shortparallel}{#2}}{#1}}
\def\O{{\mathcal O}}
\def\L{{\cal L}}
\def\M{{\cal M}}
\def\N{{\cal N}}
\def\P#1{{\mathbb P}^#1}
\def\PP{{\mathbb P}}
\def\C{{\mathbb C}}
\def\Z{{\mathbb Z}}
\def\R{{\mathbb R}}
\def\Q{{\mathbb Q}}
\def\sym{{\mathrm{Sym}}}
\setcounter{section}{-1}
\newcommand{\defref}[1]{Def.~\ref{#1}}
\newcommand{\propref}[1]{Prop.~\ref{#1}}
\newcommand{\thmref}[1]{Thm.~\ref{#1}}
\newcommand{\lemref}[1]{Lemma~\ref{#1}}
\newcommand{\corref}[1]{Cor.~\ref{#1}}
\newcommand{\exref}[1]{Example~\ref{#1}}
\newcommand{\secref}[1]{Section~\ref{#1}}

\newcommand{\qedd}{\hfill\framebox[2mm]{\ }\medskip}
\newcommand{\codim}{\textrm{codim}}
\def\proof{{\it Proof:\  \ }}

\author{Giorgio Ottaviani}
\thanks{The author is member of Italian GNSAGA-INDAM. Partially supported by the
PRIN project ``Multilinear Algebraic Geometry'' of MUR.}
\address{Department of Mathematics and Computer Science ``Ulisse Dini'', University of Florence, Italy}
\email{giorgio.ottaviani@unifi.it}

\subjclass[2020]{14-02, 14F17, 14F08, 14M15, 14M17}

\date{}
\title[Vector bundles without intermediate cohomology and the trichotomy]{Vector bundles without intermediate cohomology and the trichotomy result}
\maketitle

\begin{small}\emph{  Dedicated to Enrique Arrondo on the occasion of his
60-th birthday.}\end{small}

\begin{abstract} 
Horrocks proved in 1964 \cite{Hor} that vector bundles on $\PP^n$ without intermediate cohomology split as direct sum of line bundles. This result has been the starting point of a great research activity on other varieties, showing interesting connections with derived categories and other areas. We follow some paths into this fascinating story, which has classical roots. The story has a culmination with the trichotomy result (finite/tame/wild) for arithmetically Cohen Macaulay (ACM) varieties obtained by Faenzi and Pons-Llopis in 2021. This is an expanded version of the talk given at the conference "Homemade Algebraic Geometry" in July 2023 at Alcal\'a de Henares celebrating Enrique Arrondo's 60th birthday.

\end{abstract}
\section{Introduction}
A vector bundle $E$ on a variety $X$ polarized by $\O(1)$ has no intermediate cohomology if $H^i(X,E\otimes\O(t))=0$ $\forall t\in\Z$, $0<i<\dim X$. This brief survey about bundles without intermediate cohomology is divided into three sections. 

In the first one we recall cohomological criteria for a vector bundle on projective space $\PP^n$ and on smooth quadrics $Q_n$ to split as a sum of line bundles. Bundles on $\PP^n$
or $Q_n$ without intermediate cohomology are finitely many, up to twist with line bundles.
This means that $\PP^n$ and $Q_n$ are {\it finite}.

In the second section we consider the approach to these topics with derived categories.
The analogous facts are the Beilinson description od the derived category of $\PP^n$ and the Kapranov description of the derived category of $\PP^n$.
The description of the  derived category of Grassmannians by Kapranov  is a breakthrough in the story. Bundles on Grassmannians without intermediate cohomology make families of
arbitrarily large dimension, this means that any Grassmannian different from $\PP^n$ or $Gr(\PP^1,\PP^3)$  is {\it wild}.

A variety polarized with a very ample line bundle $\O(1)$ is called arithmetically Cohen-Macaulay (ACM) when $\O(1)$ has no intermediate cohomology 
In the third section we expose the classification of bundles on ACM varieties into the three classes finite/tame/wild, the {\it trichotomy}.

The style of the exposition is informal and concerns the development of main ideas through examples.
There are three open problems left respectively as questions \ref{quest1},  \ref{quest2},  \ref{quest3}. I thank Daniele Faenzi for several useful discussions.

\section{The Horrocks criterion, the spinor bundles, classical roots on ACM bundles}
\subsection{The Horrocks criterion on $\PP^n$ and its classical roots}

The starting point of this story is the remarkable
\begin{thm0}{Horrocks criterion (1964)}\label{thm:Hor}\cite{Hor}
Let $E$ be a vector bundle on $\PP^n(\C)$.

$E=\oplus_{i=1}^r\O(a_i)$ for some integers $a_i$ $\Longleftrightarrow$ $H^i(E(t))=0,\quad 1\le i\le n-1, \forall t\in\Z$.
\end{thm0}
\vskip 0.2cm

This result characterizes bundles on $\PP^n$ without intermediate cohomology
, called also arithmetically Cohen-Macaulay(ACM) bundles.
\vskip 0.6cm

 Note for $n=1$ the condition in the criterion is empty and we get

\begin{coro0}(Segre-Grothendieck)\label{cor:SG} Any vector bundles $E$ on $\PP^1$ split as the sum of line bundles,
that is
$$E=\oplus_{i=1}^r\O_{\PP^1}(a_i)$$
for some integers $a_i$.
\end{coro0}

Corollary \ref{cor:SG} is commonly attributed to Grothendieck although, as often happens in Mathematics, Grothendieck
was not the first to prove it.
A short history is contained in \cite[\S 2.4]{OSS}.

We list of some proofs in reverse chronological order, revising slightly the list in \cite[\S 2.4]{OSS},
where most of the original references are quoted.
\begin{itemize}
 \item Grauert-Remmert (1975), the most popular proof, available also in \cite[Theorem 2.1.1]{OSS},
\item Grothendieck (1956),
\item Hilbert (1906), in the language of matrices on $\C^*$,
\item Bellatalla (1901), his diploma thesis with C. Segre supervision \cite{Bel},
\item C. Segre, for rank $2$ (1883), a geometric argument \cite{S1}, 
\item Dedekind-Weber (1882), in algebraic setting.
\end{itemize}
\vskip 0.4cm

{For the convenience of the reader we recall the sketch of the Grauert-Remmert proof of Segre-Grothendieck Cor. \ref{cor:SG}.}
Let $E$ be a vector bundle on $\PP^1$.

$\bullet$ {\bf Step 1}, up to twist by a line bundle, we may assume $h^0(E)\neq 0$, $h^0(E(-1))=0$.
Let $\O\rig{s}E$ be a section, it has no zeroes, otherwise if $s(x)=0$ then $s$ is a section of $E\otimes I_x=E(-1)$.

$\bullet$ {\bf Step 2}, we get the exact sequence
\begin{equation}\label{eq:EF}0\rig{}\O\rig{s}E\rig{}F\rig{}0\end{equation}
where $F=\oplus_{i=1}^{r-1}\O(b_i)$ splits by induction on the rank.

$\bullet$ {\bf Step 3}, tensoring (\ref{eq:EF})  by $\O(-1)$ we get
$$\stackbelow{H^0(E(-1))}{0}\rig{}H^0(\oplus_{i=1}^{r-1}\O(b_i-1))\rig{}\stackbelow{H^1(\O(-1))}{0}$$
hence $b_i\le 0$.

$\bullet$ {\bf Step 4}, the sequence (\ref{eq:EF}) splits since $Ext^1(F,\O)=H^1(F^\vee)=0$.

For a modern exposition of Segre paper\cite{S1} we refer to \cite[\S 2, 3, 4]{GO}.
Segre worked in the setting of ruled varieties. In modern notations,  for a $2$-bundle $E$ on $\PP^1$ he considered the projection $\PP(E)\rig{\pi}\PP^1$ where $\PP(E)$ is a ruled surface.
The goal was to construct a section of $\pi$ such that its image make a subbundle of $E$.
Segre's proof resembles steps 1, 2 and 3 of the Grauert-Remmert proof and anticipates them. Step 4 needed 
a different approach.
\vskip 0.2cm
Barth and Hulek had the smart idea  in \cite{BH} to prove Theorem \ref{thm:Hor} by induction on $n$ by taking Corollary \ref{cor:SG}
as starting point of the induction. The bundle $E$ restricted to a line splits as $\oplus_{i=1}^r\O(a_i)$, and one lifts the isomorphism 
$E\to \oplus_{i=1}^r\O(a_i)$ from the line to the whole $\PP^n$, by the cohomological vanishing that we have in the assumptions.
This approach has been nicely exposed in \cite[Theorem 2.3.1]{OSS}.

\vskip 0.2cm

Giuseppe Gherardelli wrote in 1943 the following remarkable result,  under Severi advice.
It can be considered the ancestor of Horrocks criterion.
\begin{thm0}[G. Gherardelli]\cite{Gh}\label{thm:gh}
Let $C\subset\PP^3$ be smooth.
Then $C$ it is a complete intersection if and only if the following two conditions hold
\begin{itemize}
\item{}(i) $C$ is subcanonical, that is $K_C=\O(e)_{|C}$ for some $e\in\Z$
\item{}(ii) $C$ is projectively normal, that is the restriction map

$H^0(\PP^3,\O(k))\to H^0(C,\O(k))$ is surjective.
\end{itemize}
\end{thm0} 
Actually the assumptions in \cite{Gh} on $C$ were weaker, but Gherardelli can conclude that $C$ is smooth. In the setting of ACM bundles, the assumption that makes Theorem \ref{thm:gh} to work is that $C$ is a
locally CM curve .
After Serre correspondence (\cite[chap. 1 \S 5]{OSS} and \cite[Theorem 8.21A (e)]{Hart} ), any locally CM subcanonical curve in $\PP^3$ is obtained as the zero locus of a section
of a rank $2$ bundle $E$, which splits as the sum of two line bundles if and only if $C$ is a complete intersection.
In this optic, Gherardelli Theorem for curves in $\PP^3$ is equivalent to Horrocks criterion for $2$-bundles on $\PP^3$, indeed the condition (ii) is equivalent to the vanishing
 $H^1(I_{C,\PP^3}(k))=0$ which from the sequence
$$0\rig{}\O_{\PP^3}\rig{}E\rig{} I_{C,\PP^3}(c_1(E))\rig{}0$$ 
is equivalent to $H^1(E(*))=0$, and the second vanishing $H^2(E(*))=0$ follows from Serre duality.

\subsection{Spinor bundles and the splitting criterion on quadrics}\label{subsec:spinor}
On {\it quadrics} there is a notable class of bundles without intermediate cohomology beyond line bundles, the {\it spinor bundles}.

The linear spaces contained in a smooth quadric $Q_n\subset\PP^{n+1}$  have maximal dimension $\lfloor n/2\rfloor$.
Linear spaces contained in a quadric are called isotropic, since a linear space $\PP(L)$ is contained in $Q_n$ if and only if
$Q(v,w)=0$ $\forall v, w\in L$ where $Q(-,-)$ is the bilinear form associated to the quadric.

 Let for simplicity $n=2k+1$ ($n$ odd). The spinor variety $S_{k+1}$ parametrizes the isotropic linear spaces of maximal dimension
 \begin{equation}\label{eq:spinor}
S_{k+1}=\left\{\PP^k\in Gr(\PP^k,\PP^{2k+2})|\PP^k\subset Q_{2k+1}\right\}
\end{equation}
and it has easily seen to have dimension $\frac{(k+1)(k+2)}{2}$. The first cases are $S_1=\PP^1$ and $S_2=\PP^3$. By the Pl\"ucker embedding, $S_{k+1}$ is embedded in the large projective space $\PP^N$ with $N={{2k+3}\choose {k+1}}-1$.   The striking fact is that ${\mathrm Pic}(S_{k+1})=\Z$ and
$$\O_{Gr}(1)_{|S_{k+1}}=\O_{S_{k+1}}(2)$$

so that the spinor variety $S_{k+1}$ has another embedding corresponding to the ``half'' line bundle $\O(1))$ (replacing $\O(2)$) in the projective space of a vector space of smaller dimension $2^{k+1}$, the famous spin representation
$(\mathrm{Spin}\C^{2k+1})$, that we recall briefly in Remark \ref{rem:spin_nutshell}.

We have the incidence variety
$$\mathbb F\colon =\left\{(x,\PP^k)\in Q_{2k+1}\times S_{k+1}| x\in\PP^{k}\right\}$$
with the two projections $p_1$ and $p_2$ and
$S^\vee=(p_1)_*(p_2)^*\O(1)$ is the dual of the spinor bundle. Since the tangent space $T_xQ_{2k+1}$ meets
$Q_{2k+1}$ in a cone with vertex $x$ over a smaller quadric $Q_{2k-1}$, the fibers of $p_1$ are isomorphic to
a smaller spinor variety $S_k$. This proves that ${\mathrm rank}S=h^0(S_k,\O(1))=2^k$. 

 On even dimensional quadrics $Q_{2k}$ there are two families of maximal isotropic subspaces $\PP^k$, both parametrized by the same spinor variety $S_k$ encountered in the odd case,
  and accordingly two spinor bundles $S'$, $S''$ of rank $2^{k-1}$. In \cite{Ott88} it is described
how the spinor bundles restrict to hyperplane quadrics, this argument allows to prove directly
by induction that ${\mathrm rank}S=2^k$.

It is useful to compare the following facts.
\begin{itemize}
\item{}
The Grassmannian $Gr(\PP^k, \PP^n)$ is the unique closed $SL(n+1)$-orbit in $\PP(\wedge^{k+1}\C^{n+1})$.
\item{}
 The Veronese variety $v_d(\PP^{n})$ is the unique closed $SL(n+1)$-orbit in $\PP(\mathrm{Sym}^{d}\C^{n+1})$.
\item{}
 The spinor variety $S_k$ is the unique closed $Spin(2k+1)$-orbit in $\PP(\mathrm{Spin}\C^{2k+1})$.
 \end{itemize}

\begin{rem0}\label{rem:spin_nutshell}
We take the opportunity to expose the spin representation in a nutshell,
as can be taught at undergraduate level.
In an algebra course, $\C$ is commonly introduced as $\R\oplus\R i$ where you {\it declare} $i^2=-1$.

 In a linear algebra course, $\C$ may be introduced as the $\R$-algebra 
 $$\C=\left\{\begin{bmatrix}a&-b\\b&a\end{bmatrix}| a, b\in\R\right\}$$, where $i=\begin{bmatrix}0&-1\\1&0\end{bmatrix}$ {\it satisfies} $i^2=-1$.
This has the advantage that the imaginary  unit  is a concrete real $2\times 2$ matrix.

 In a second course, you may introduce the quaternion as the $\C$-algebra
$${\mathbb H}=\left\{\begin{bmatrix}z&-w\\\overline{w}&\overline{z}\end{bmatrix}| z, w\in\C\right\}$$
and now you pose
$${\mathbf i}=\begin{bmatrix}i&0\\0&-i\end{bmatrix}, {\mathbf j}=\begin{bmatrix}0&-1\\1&0\end{bmatrix},
{\mathbf  k}=\begin{bmatrix}0&-i\\i&0\end{bmatrix}$$
 they satisfy the usual rules, and they are a realization of the Pauli matrices (with a different normalization).
The quaternion conjugate corresponds to the operation ${\overline{A}}^t$ on the  matrix $A$
and this allows to define quaternion norm as $|q|^2=q\overline{q}$ for $q\in{\mathbb H}$.

  Now ${\mathbb H}=\R\oplus{\mathbb I}$ contains $S^3=\{q| |q|=1\}$ which shows the remarkable fact that $S^3$ has a group structure, indeed it is a real Spin group.
 
  Construct the map
 $$\begin{array}{ccc}S^3&\to& SO(3)=SO({\mathbb I})\\
 q&\mapsto&(x\mapsto q^{-1}xq)\end{array}$$
 which is $2:1$ group homomorphism, since the fibers consist of $\{-q, q\}$.
 \vskip 0.6cm
 
It follows that $S^3$ is the universal covering of $SO(3)$, we pose $S^3=Spin(3)$ and the natural inclusion
 $S^3\subset GL(2,\C)$ as quaternions of norm $1$ identifies $S^3=SU(2)$ and {\it is} the $2$-dimensional
 real spin representation of $Spin(3)$. On complex numbers we have, analogously, $Spin(3)=SL(2)$.

The general $Spin$ representation follows this path by replacing the algebra of quaternions ${\mathbb H}$ with the more general Clifford algebra.
\end{rem0}

The spinor bundle is ACM, so that to get a splitting criterion analog to Horrocks criterion, it is necessary to add a further condition,
as in the following.

\begin{thm0}{{Splitting Criterion on Quadrics, \cite{Ott89}}
{Let $E$ be a vector bundle on $Q_n$, let $S$ be a spinor bundle. Then

$E=\oplus_{i=1}^r\O(a_i)$ $\Longleftrightarrow$ $\left\{\begin{array}{lr}H^i(E(*))=0,& 1\le i\le n-1\\
H^{n-1}(E\otimes S(*))=0,&\end{array}\right.$}}
\end{thm0}

Nevertheless, Horrocks criterion can be generalized in a second way, by describing the class of ACM bundles on quadrics.

\begin{thm0} {{Characterization of ACM bundles on Quadrics, Kn\"orrer, \cite{Kn}}
{Let $E$ be an indecomposable vector bundle on $Q_n$. Then
\vskip 0.2cm

$E=\O(a)$ or $E=S(b)$  $\Longleftrightarrow$ $H^i(E(*))=0, 1\le i\le n-1$}}
\end{thm0}

 It is easy to check  that
$$\textrm{Splitting Criterion on\ }Q_n\Longleftarrow \textrm{Kn\"orrer Characterization}$$

 Enrique patiently explained to me in 1989, at the conference "Projective varieties" in Trieste, that also the converse implication holds !

His technique, by successive ``killing of $H^1$'', will be exposed later in \cite{AG} that we recall in a while, see Theorem \ref{thm:AG}.

\subsection{Improvements of Horrocks criterion for small rank}
The results in this section improve the Horrocks criterion when the rank of $E$ is relatively small.

\begin{thm0}[Evans-Griffith criterion (1981)]\label{thm:EG}\cite{EG}
{Let $E$ be a vector bundle on $\PP^n(\C)$ of rank $r\le n$.

$E=\oplus_{i=1}^r\O(a_i)$ $\Longleftrightarrow$ $H^i(E(*))=0,\quad 1\le i\le  r-1$}
\end{thm0}

Note for $2$-bundles only the vanishing of $H^1(E(*))$ characterizes the splitting of $E$. A simple proof of Theorem \ref{thm:EG} as a consequence of Le Potier Vanishing Theorem 
has been given by L. Ein in \cite{E}.

\vskip 0.5cm
 \begin{thm0}[Kumar-Peterson-Rao (2003)]\label{thm:KPR}\cite{KPR}
Let $\mathrm{rk}(E)\le n-2$, $n\ge 3$.
Then $E$ splits if and only if $H^i(E(*))=0\quad 2\le i\le n-2$.
\end{thm0}

The cohomological bound in Theorem \ref{thm:KPR} has been called {\it  inner cohomology}
so that this Theorem characterizes bundles without { inner cohomology}.
Note in particular a $2$-bundle $E$ on $\PP^4$ splits if and only if $H^2(E(*))=0$.

The following result generalizes Theorem \ref{thm:KPR} to quadrics.
\begin{thm0}[Ancona-Peternell-Wisniewski 1994, Malaspina 2009]\cite{APW, Mal09}
Let $E$ be a rank $2$ bundle on $X=\PP^n$ or $X=Q_n$, $n\ge 4$. TFAE
\begin{enumerate}
\item $H^i(E(*))=0$ for $2\le i\le n-2$ (without inner cohomology)
\item $\PP(E)$ is Fano
\item $E=\O(a)\oplus\O(b)$, $E$ is a spinor bundle on $Q_4$, and two further examples.
\end{enumerate}
\end{thm0}
\vskip 0.4cm

 \begin{quest0}\label{quest1} Are the first two items equivalent on other Fano varieties of dimension $\ge 4$ ? Note that the proof of the equivalence on $\PP^n$ and $Q_n$ is indirect, through the third item.
\end{quest0}

Enrique Arrondo and Laura Costa\cite{AC} classified rank-2 indecomposable ACM bundles over Fano
$3$-folds of index 2 by Hartshorne–Serre correspondence. 
\vskip 0.2cm

\begin{thm0}[Arrondo-Costa]\cite{AC}\label{thm:AC} There
are three families of rank $2$ indecomposable ACM bundles on the smooth Fano $3$-fold of index $2$ ($V_3$, $V_4$ and $V_5$), having a section of a suitable twist vanishing respectively on
\begin{enumerate}
\item a line
\item a conic 
\item a projectively normal elliptic curve (in the threefold).
\end{enumerate} 
\end{thm0}

This result has been generalized to very Fano $3$-fold by Brambilla and Faenzi in \cite{BF}
and to $V_{22}$ by Arrondo and Faenzi in \cite{AF},

\subsection{Grassmannians and beyond, the appearance of many ACM bundles}\label{subsec:gras}

It remained an open problem to characterize ACM bundles on Grassmannians.
Since the universal and the quotient bundle are ACM one could at first hope that they are the only ones or just a few ones.
The following result by Arrondo, Graña in 1999 destroyed this hope.

Arrondo - Grana (1999) prove that ACM bundles on $Gr(\PP^1,\PP^4)$ make big families,
in particular there exist {\it non homogeneous} ACM bundles. 
\vskip 0.2cm
\begin{thm0}[Arrondo, Graña]\cite{AG}\label{thm:AG}
There are families of positive dimension of vector bundles $E$ on ACM bundles on $Gr(\PP^1,\PP^4)$ .
Some of these examples are non homogeneous under $SL(5)$-action.
\end{thm0}

This result is a turning point in the story since shows that the structure of ACM bundles can be in general
quite complicated.
It is worth that the class of ACM bundles on $Gr(\PP^1,\PP^4)$ is explicitly described, up to extensions of known bundles.

What is in embryo in \cite{AG} is the concept that $Gr(\PP^1,\PP^4)$ is wild, that is, it carries families of ACM bundles with arbitrarily large dimension. This was understood later and popularized by Hartshorne, see the section \ref{sec:trichotomy}.
Wilderness of any Grassmannian different from projective spaces and $Gr(\PP^1,\PP^3)$ was definively settled by Costa and Mirò-Roig in \cite{CMR16}.

\begin{rem0} I recall several discussions with Enrique Arrondo at the time of \cite{AG}, when we did not yet image the wilderness of most Grassmannians, and we were already impressed by the non homogeneous examples.
\end{rem0}

I have started this subsection with  the results of \cite{AG}, since I wished to show immediately the appearance of many ACM bundles. The chronological order is a bit different, since 10 years before it was proved the following splitting criterion.
Using the language of Schur functors, for a bundle $E$ and a Young diagram $\alpha$ it is defined the bundle $S^{\alpha}E$.
When $\alpha$ is  a horizontal Young diagram like $\alpha=\yng(4)$ then $S^\alpha E={\mathrm Sym}^4E$,
while if $\alpha$ is  a vertical Young diagram like $\alpha=\yng(1,1,1)$ then $S^\alpha E=\wedge^3E$.
General Young diagrams act by ``mixing'' the symmetric and skew-symmetric actions.

It turns out, by using Bott Theorem, that the universal bundle $U$ is ACM and $S^\alpha U$ is ACM
if and only if $\alpha$ has at most $n-k-1$ columns.

\begin{thm0}[O, 1989]\label{thm:ottgras}
{Let $E$ be a vector bundle on $G=Gr(\PP^k, \PP^n)$.

$E=\oplus_{i=1}^r\O(a_i)$ $\Longleftrightarrow$ $H^i(E(*)\otimes S^\alpha U)=0,\quad 1\le i\le \dim G-1$ 

where $\alpha$ is a Young diagram with at most $n-k-1$ columns, 
${U}$ is the universal bundle of rank $k+1$.}
\end{thm0}

The proof is by induction on $k$, in Grauert-Remmert style.

Arrondo-Malaspina \cite{AM} and Arrondo-Tocino \cite{AT}
find a characterization of symmetric powers of quotient bundle on $Gr(\PP^1,\PP^n)$.
We state explicitly a result from the second paper.

\begin{thm0}[Arrondo-Tocino]\cite{AT}
Let $E$ be a bundle on $G=Gr(\PP^1, \PP^n)$.

$E=\oplus_{i=1}^r\O(a_i)$ $\Longleftrightarrow$ $\left\{\begin{array}{cc}H^i(E\otimes\mathrm{Sym}^iQ(*))=0&1\le i\le n-2\\
H^i(E\otimes\mathrm{Sym}^{2n-3-i}Q(*))=0&n-1\le i\le 2n-3\end{array}\right.$
\end{thm0}

The main result of \cite{AT} is the determination of a finite subset $S_k$ of pairs such that, given $0\le k\le n-2$,
$E$ is a direct sum of twists of $\O, Q, \ldots, \mathrm{Sym}^kQ$ if and only if  $H^j(E\otimes\mathrm{Sym}^iQ(*))=0$ for $(i,j)\in S_k$.
\vskip 0.8cm

{We end this section with the open case of Lagrangian Grassmannians}
Let $\C^{2k+2}$ be equipped with a nondegenerate skew form $J$. Let $LG(k)$ be the Lagrangian Grassmannian of 
maximal isotropic $\PP^k\subset \PP^{2k+1}$, that is
\begin{equation}\label{eq:LGk}LG(k)=\left\{L\subset \C^{2k+2}| \dim L=k+1, J(v,w)=0\quad\forall v, w\in L\right\}.
\end{equation}
Let $E$ be a vector bundle on $LG(k)$.
A naive attempt to generalize the splitting criterion on Grassmannian to $LG(k)$  is the following.

$E=\oplus_{i=1}^r\O(a_i)$ $\overbrace{\Longleftrightarrow}^{?}$ $H^i(E(*)\otimes S^\alpha U)=0, 1\le i\le dim LG(k)-1$
where $\alpha$ has at most $k$ columns.
\vskip 0.4cm

 While $\Longleftarrow$ is true, the implication $\Longrightarrow$  is true on $LG(1)=Q_3$ but it is false for $k=2$.

 Oeding and Macias Marques have made, during 2009 Pragmatic, an attempt which works, but only partially. 

\begin{thm0}[Oeding-Macias Marques]\cite{OM}
Let $E$ be a vector bundle on $LG(k)$.
Let $C_k(E)$ be the vanishing condition
$$H^i(LG(k),\wedge^{j_1}U\otimes\ldots\wedge^{j_k}U\otimes E(*))=0, 0\le j_q\le q+1, i=\sum_{q=1}^kj_q$$
\begin{enumerate}
 \item If $C_k(E)$ holds then $E$ splits as a sum of line bundles.
 \item If $E$ splits as a sum of line bundles and $k\le 6$ then $C_k(E)$ holds.
 \item Let $k=7$ and $E=\O$. Then $C_7(E)$ does NOT hold.
\end{enumerate}
\end{thm0}
\begin{quest0}\label{quest2}
 Can be found a cohomological splitting criterion on $LG(k)$ for any $k$ ? One should likely use irreducible
homogeneous bundles (with the language of maximal weights) and not just products of wedge powers.
\end{quest0}

\section{The approach with the  derived category}

\subsection{The Grothendieck groups and a bit of K-theory}

Let $X$ be a smooth variety. Let $K^0(X)$ (resp. $K_0(X)$) be the free group generated by isomorphism classes of vector bundles (resp. coherent $\O_X$-modules) modulo the relation
$[A]-[B]+[C]=0$
when we have the exact sequence
$$0\rig{}A\rig{}B\rig{}C\rig{}0$$.

$K^0(X)$ is called the Grothendieck group of vector bundles on $X$ and $K_0(X)$ is called the Grothendieck group of coherent sheaves on $X$.

The following is one of the main results by Grothendieck in his  SGA6, we recommend the readable proof by Murre in \cite[\S 2.10]{Mur}.
\begin{thm0}[Grothendieck]

\begin{enumerate}
\item $K^0(X) = K_0(X)$
\item the Chern character gives a ring isomorphism
$$ch\colon K^0(X)\otimes\Q\to CH(X)_{\Q}$$
where $CH(X)_{\Q}$ is the rational Chow ring.
\end{enumerate}
\end{thm0}

For any rational homogeneous variety $X=G/P$, every element of the Lie algebra ${\mathfrak g}$ gives an element of the tangent space at $[P]$
which, with $G$-action, defines a section of $TX$. It is a well known fact that all sections are defined in this way,
that is $H^0(TX)={\mathfrak g}$. 

This fact is important since it allows to compute the zero locus of a general section of $TX$, hence, by Gauss-Bonnet Theorem,
the topological Euler characteristic $\chi(X)$. It is well known that the category of $G$-bundles on $G/P$ is equivalent to the category of $P$-modules. The tangent bundle $T (G/P)$ is associated to an {\it irreducible} representation of $P$ 
just in a few cases which are called irreducible Hermitian symmetric. They are the three classical cases listed in Proposition \ref{prop:chi}, smooth quadrics (where $\chi Q_n$ is equal to $n+1$ if $n$ is odd and to $n+2$ if $n$ is even)
and just two exceptional cases for the groups $E_6$ and $E_7$ that we skip for simplicity. The following result is well known (indeed all Betti numbers are known),
but we think it is worth to expose the simple and unified proof for these three cases.

\begin{prop0}\label{prop:chi}
\begin{enumerate}
\item{}$\mathrm{rank} K_0\left(Gr(k,n)\right) = \chi\left(Gr(k,n)\right)={n\choose k}$
\item{}$\mathrm{rank} K_0\left(LG(n+1)\right) = \chi LG(n+1)=2^n \left( \textrm{see} (\ref{eq:LGk}) \right)$
\item{}$\mathrm{rank} K_0\left(Spin_n\right) = \chi Spin_n=2^{n}  \left( \textrm{see} (\ref{eq:spinor}) \right)$
\end{enumerate}
\end{prop0}
\begin{proof}
In all three cases a section is given by a map $A\colon V\to V$.
In case (i) it is general, in case (ii) it preserves the quadric $Q$ (that is $AQ+QA^t=0$), in case (iii)
it preserves the skew-symmetric form  $J$ (that is $AJ+JA^t=0$).
We may assume, by genericity, that $A$ has $n$ distinct eigenvectors with nonzero eigenvalues.
 In case (i), the section $A$ vanishes on the linear spaces $\C^k$
which are span of $k$ among the $n$ eigenvalues, they are ${n\choose k}$.
In case (ii), we have the further restriction that the subspace is isotropic.
The eigenvectors of $A$ satisfy $vQv^t=0$. 
We may assume we have an orthonormal basis $\left\{ e_1,\ldots, e_n, f_1,\ldots, f_n\right\}$
with $Q(e_i,f_j)=\delta_{ij}$, $Q(e_i, e_j)=Q(f_i, f_j)=0$. 
The isotropic subspaces are generated by choosing exactly one between $e_i$ and $f_i$, for $i=1,\ldots n$,
so we have $n$ binary choices and $2^n$ subspaces.
Hence the isotropic $C^n$ are $2^n$ ($2^{n-1}$ in one family) and the isotropic $C^k$ are ${n\choose k}2^k$.
In case (iii) the construction is analogous. Now all vectors are isotropic,
in particular the eigenvectors are isotropic.
\end{proof}

\subsection{The derived category of $\PP^n$ and $Q_n$}
We refer to \cite{GM} for basics on derived categories. See also \cite{AO89} for a short introduction.
The derived category $D_b(X)$ of bounded complexes of sheaves on $X$ can be thought as an {\it enhancement}
of $K_0(X)$. Recall that a sequence of bundles $(E_0,\ldots, E_n)$ on a variety $X$ is called exceptional if $Hom(E_i,E_i)=\C$
and $Ext^m(E_i, E_j)=0$ for $i\ge j$, $m>0$ or $i>j$, $m\ge 0$ and it is called full if it generates  $D_b(X)$.

\begin{thm0}[Beilinson, 1978]\cite{Bei, AO89}
The derived category $D_b(\PP^n)$
has a full  exceptional sequence (can be understood as a basis) $\langle \O(-n),\ldots, \O\rangle$
\end{thm0}

\begin{coro0}
$K_0(\PP^n)$ is freely generated by the classes 
$$[\O(-n)],\ldots, [\O]$$
\end{coro0}

Note the Chow ring of $\PP^6$ has the following graphical representation
\begin{center}
\setlength{\unitlength}{3.5mm}
\begin{picture}(30,9)(-8,2)
\multiput(-4.2,7.7)(2,0){7}{$\bullet$}
\put(-4.2,6.8){$0$}
\put(-2.2,6.8){$1$}
\put(-0.2,6.8){$2$}
\put(1.8,6.8){$3$}
\put(3.8,6.8){$4$}
\put(5.8,6.8){$5$}
\put(7.8,6.8){$6$}
\put(-4,8){\line(1,0){12}}
\end{picture}
\end{center}
Beilinson in \cite{Bei} considers the Segre product $\PP(V)\times\PP(V)$ with the two projections on the factors $p_1$, $p_2$.
Let 
\begin{equation}\label{eq:BeilE}E=p_1^*Q\otimes p_2^*\O(1)\end{equation}
where $Q$ is the quotient bundle.
We have the identification $H^0(E)=End(V)$ and the section $s$ corresponding to the identity endomorphism vanishes on the diagonal 
$\Delta\subset \PP(V)\times\PP(V)$. The sheaf  $\O_{\Delta}$ is equivalent in the derived category to its resolution given by the Koszul complex of $s$.
Beilinson smart idea is to compare $E={Rp_1}_*(p_2^*E_{|\Delta}$
with the derived functor $R{p_1}_*$ applied to the Koszul complex of $s$ tensored by $p_2^*E$.

Ten years later, Kapranov extends this construction to quadrics and Grassmannians in \cite{Kap}.
Let state first his main result for quadrics.

\begin{thm0}[Kapranov, 1988]\cite{Kap}
The derived category $D_b(Q_n)$
has a full  exceptional sequence $\langle \O(-n+1),\ldots, \O, S\rangle$ for $n$ odd
and by
$\langle \O(-n+1),\ldots, \O, S', S''\rangle$ for $n$ even,
where $S$, $S'$, $S''$ are the spinor bundles introduced in \S\ref{subsec:spinor}.
\end{thm0}

 \begin{coro0}
$K_0(Q_n)$ is freely generated by the classes 
$$[\O(-n+1)],\ldots, [\O], [S]$$ for $n$ odd, by the classes 
$$[\O(-n+1)],\ldots, [\O], [S'], [S'']$$ for $n$ even.

\end{coro0}
Note the Chow ring of $Q_n$ is analogous to the one of $\PP^n$ when $n$ is odd, in the sense that the Betti numbers are the same, but
the powers $H^i$, with $H$ the hyperplane divisor, generate $CH^i$ for $i<\lfloor n/2\rfloor$ and are equal to twice the generator for $i\ge\lfloor n/2\rfloor$.
When $n$ is even the middle part has dimension $2$, as it is shown in the following diagram of $Q_6$.  It is interesting to remark that the splitting of the middle cohomology testifies  the
presence of the two spinor bundles, replacing the single spinor bundle for the odd case.
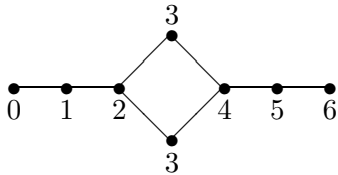
\begin{center}
\setlength{\unitlength}{3.5mm}
\begin{picture}(30,12)(-8,2)
\multiput(-4.2,7.7)(2,0){3}{$\bullet$}
\put(-4.2,6.8){$0$}
\put(-2.2,6.8){$1$}
\put(-0.2,6.8){$2$}
\put(1.8,4.8){$3$}
\put(1.8,10.4){$3$}
\put(3.8,6.8){$4$}
\put(5.8,6.8){$5$}
\put(7.8,6.8){$6$}
\put(-4,8){\line(1,0){4}}
\put(1.8,9.7){$\bullet$}
\put(1.8,5.7){$\bullet$}
\multiput(3.8,7.7)(2,0){3}{$\bullet$}
\put(4,8){\line(1,0){4}}
\put(0,8){\line(1,1){2}}
\put(0,8){\line(1,-1){2}}
\put(2,10){\line(1,-1){2}}
\put(2,6){\line(1,1){2}}
\end{picture}
\end{center}

\subsection{Cohomological characterizations via derived category}
In \cite{AO91} we found an alternative proof of Horrocks criterion (Theorem \ref{thm:Hor}) based on Beilinson description of $D_b(\PP^n)$.
The proof allows to extend the characterization to any coherent sheaf as in the following.

\begin{thm0}[Ancona-O, 1991]\cite[Theorem 2.2]{AO91}
{Let $E$ be a coherent sheaf on $\PP^n(\C)$.

$E=\oplus_{i=1}^r\O(a_i)\oplus {\mathcal S}$ $\Longleftrightarrow$ $H^i\left(E(*)\right)=0,\quad 1\le i\le  n-1$}

where ${\mathcal S}$ is supported on a $0$-dimensional scheme.
\end{thm0}

 An analogous statement\cite[Corollary 6.2]{AO91} holds for quadrics, using Kapranov description of $D_b(Q_n)$.

It is now time to state the remarkable description of the derived category of Grassmannians given by Kapranov.
Kapranov replaces the bundle $E$ in (\ref{eq:BeilE}) with the tensor product
$p_1^*Q\otimes p_2^*U^\vee$ on $Gr(\PP^k,\PP^n)\times Gr(\PP^k,\PP^n)$ where $Q$ is the quotient bundle and $U$ is the universal bundle.
Again there is a section vanishing on the diagonal and the whole construction of \cite{Bei} can be repeated at the price of more technical
stuff involving Schur functors. The final result is the following.

 \begin{thm0}[Kapranov, 1988]\label{thm:kaprgras}\cite{Kap}
The derived category $D_b(Gr(\PP^k,\PP^n)$
 has the full  exceptional sequence $\langle S^\alpha U\rangle$ for 
$\alpha$ Young diagram with at most $n-k-1$ columns, where $U$ is the universal bundle. In particular such $S^\alpha U$ generate 
$K_0(Gr(\PP^k, \PP^n)$.
\end{thm0}

The Hasse diagram of Schubert varieties in $Gr(\PP^1,\PP^4)$, where the vertices generate the Chow ring $K_0(Gr(\PP^k, \PP^n)$ is the following (see also \cite[\S 4]{OR}):

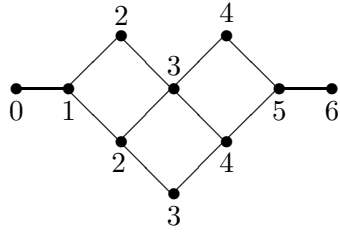
\begin{center}
\setlength{\unitlength}{3.5mm}
\begin{picture}(30,12)(-8,2)

\multiput(-2.2,7.7)(2,0){2}{$\bullet$}
\multiput(1.8,9.7)(4,0){2}{$\bullet$}
\multiput(1.8,5.7)(4,0){2}{$\bullet$}
\multiput(3.8,7.7)(4,0){2}{$\bullet$}
\multiput(9.8,7.7)(2,0){1}{$\bullet$}
\put(3.8,3.7){$\bullet$}

\put(-2.2,6.8){$0$}
\put(-0.2,6.8){$1$}
\put(1.7,4.8){$2$}
\put(1.8,10.4){$2$}
\put(3.8,8.5){$3$}
\put(3.8,2.8){$3$}
\put(5.8,4.8){$4$}
\put(5.8,10.5){$4$}
\put(7.8,6.8){$5$}
\put(9.8,6.8){$6$}

\put(-2,8){\line(1,0){2}}
\put(0,8){\line(1,1){2}}
\put(0,8){\line(1,-1){4}}
\put(2,10){\line(1,-1){4}}
\put(6,10){\line(1,-1){2}}
\put(2,6){\line(1,1){4}}
\put(4,4){\line(1,1){4}}
\put(8,8){\line(1,0){2}}

\end{picture}
\end{center}

The corresponding partitions $\alpha$ corresponding to the generators  $S^\alpha U$ of 
$K_0(Gr(\PP^k, \PP^n)$ are in the following table

$$
\begin{array}{|c|c|c|c|c|c|c|}
\hline
0&1&2&3&4&5&6\\
\hline
&&\yng(2)&&\yng(2,2)&&\\
&\yng(1)&&\yng(2,1)&&\yng(2,2,1)&\yng(2,2,2)\\
&&\yng(1,1)&&\yng(2,1,1)&&\\
&&&\yng(1,1,1)&&&\\
\hline
\end{array}$$

We warn the reader that the Hasse diagram of Schubert varieties in $Gr(\C^3,\C^6)$ is no more planar.

\def\niente2{The Hasse diagram of Schubert varieties in $S_{10}$ is the following:

\begin{center}
\setlength{\unitlength}{3.5mm}
\begin{picture}(30,12)(-8,2)

\multiput(-4.2,7.7)(2,0){3}{$\bullet$}
\multiput(1.8,9.7)(4,0){3}{$\bullet$}
\multiput(1.8,5.7)(4,0){3}{$\bullet$}
\multiput(3.8,7.7)(4,0){3}{$\bullet$}
\multiput(13.8,7.7)(2,0){2}{$\bullet$}
\put(7.8,11.7){$\bullet$}
\put(3.8,3.7){$\bullet$}

\multiput(-4.2,6.8)(2,0){3}{$1$}
\put(1.7,4.8){$1$}
\put(1.8,10.4){$1$}
\put(3.8,8.5){$2$}
\put(3.8,2.8){$1$}
\put(5.8,4.8){$3$}
\put(5.8,10.5){$2$}
\put(7.8,6.8){$5$}
\put(7.8,12.5){$2$}
\put(9.8,4.8){$5$}
\put(9.8,10.5){$7$}
\multiput(11.7,6.8)(2,0){3}{$12$}

\put(-4,8){\line(1,0){4}}
\put(0,8){\line(1,1){2}}
\put(0,8){\line(1,-1){4}}
\put(2,10){\line(1,-1){4}}
\put(6,10){\line(1,-1){4}}
\put(8,12){\line(1,-1){4}}
\put(2,6){\line(1,1){6}}
\put(4,4){\line(1,1){6}}
\put(10,6){\line(1,1){2}}
\put(12,8){\line(1,0){4}}

\end{picture}
\end{center}

Note there are exactly $16$ vertices.
}

There is an evident analogy between the generators of $D_b(Gr(\PP^k,\PP^n)$
in Theorem \ref{thm:kaprgras} and the splitting criterion on Grassmannians in Theorem \ref{thm:ottgras}.

\begin{quest0}\label{quest3} Is there a proof of the splitting criterion Theorem \ref{thm:ottgras} by using the
derived category approach, analog to \cite{AO91} for the case of projective spaces and quadrics ?
Some results for Segre products are in \cite{Mal08, CMR05}.
\end{quest0}

\begin{rem0} After the work of Kapranov on Grassmannians, there is a natural folklore conjecture (see \cite{Fon})
that the bounded derived category
of a rational homogeneous variety $G/P$ admits a full exceptional collection of $G$-homogeneous bundles. 
\end{rem0}

The conjecture of previous remark is open even for irreducible Hermitian symmetric spaces, but the case of isotropic Lagrangian
Grassmannian $LG(n)$ has been recently solved.  Kuznetsov and  Polishchuk 
constructed collections on $LG(n+1)$ of the right length $2^n$ in \cite{KP}. 

Fonarev in \cite{Fon} shows the fullness of the exceptional collections constructed by Kuznetsov
and Polishchuk and gives a geometric description. The full  collection consists of $G$-equivariant vector bundles.

Beyond the rational homogeneous varieties $G/P$, we conclude this section with some of the many results on the Fano case.
Recall the Fano $3$-fold $V_5$ ( already met exposing \cite{AC} in Theorem \ref{thm:AC}) is $SL(2)$-quasi homogeneous. The derived category of $V_5$ has been found by Orlov in \cite{Orl91}. Faenzi in \cite{Fa05}
found more mutations and applied them to classify the ACM bundles.
 \begin{thm0}[Faenzi, Orlov]
The Fano $3$-fold of index $2$ $V_5$ has the full  exceptional sequence
$$\langle \O(-1), U, Q^\vee, \O\rangle$$
where $U$ is the rank $2$ universal bundle and $Q$ is the rank $3$ quotient bundle,
both $SL(2)$-invariant.
\end{thm0}
\vskip 0.4cm

 Kuznetsov in \cite{Kuz96} have found  a similar result for the Fano $3$-fold $V_{22}$, which has index $1$,
but it has a member of the deformation class which is $SL(2)$-quasi homogeneous. Faenzi in \cite{Fa07}
found again more mutations and applied them to classify the ACM bundles.

Kuznetsov found also the decomposition of the {derived category of cubic $3$-fold}.

\begin{thm0}[Bondal-Orlov, Kuznetsov]\cite{Kuz04}
Let $X\subset\PP^{n+1}$ be a smooth cubic hypersurface.
$D_b(X)$ has a semiorthogonal decomposition
$$\langle {\bf T}_X,\O, \O(1),\ldots, \O(n-2)\rangle$$
\end{thm0}

 \begin{thm0}[Lahoz-Macr\`\i-Stellari, 2015]\cite{LMS}
Let $X\subset\PP^4$ be a smooth cubic $3$-fold.
Any ACM stable bundle on $X$ has a twist which belongs to ${\bf T}_X$.
\end{thm0}
This result is applied to show non-emptiness and smoothness of moduli space of Ulrich bundles 
(of rank $\ge 2$) on a smooth cubic $3$-fold.

\section{The trichotomy  for ACM varieties}\label{sec:trichotomy}
Buchweitz, Greuel and Schreyer proved in \cite[Theorem C]{BGS} that hypersurfaces in $\PP^n$, $n\ge 2$ of degree $d\ge 3$ have infinite non isomorphic ACM bundles.
This is the first step toward the wilderness of such hypersurfaces as defined below.
I heard the first suggestions of the trichotomy finite/tame/wild in a memorable talk by Hartshorne in Turin in 2005. In that conference Hartshorne suggested to study
ACM bundles of rank two on cubic surfaces, which was carried out by Faenzi in \cite{Fa08}.
Then Casanellas and Hartshorne proved that cubic surfaces are wild in \cite{CH}, where they introduced the terminology finite/tame/wild for ACM sheaves.
It was discovered already in \cite{BGS}  that the
quadric cone in $\PP^3$ has a countable discrete set of indecomposable ACM bundles
of rank $2$, so making the tame case for singular varieties quite subtle.
\begin{thm0}[Faenzi - Pons Llopis, 2019, Trichotomy \cite{FP}.]
Any reduced ACM closed subscheme $X\subset\PP^n$ of positive dimension falls in
exactly one of the following classes:
\begin{itemize}
 \item{\bf Finite:} there are only finitely many indecomposable ACM sheaves over $X$ up to isomorphism.
 \item{\bf Tame:} for any given rank $r$, the moduli space of indecomposable non-isomorphic ACM
sheaves of rank $r$ is a finite or countable union of points or curves.
 \item{\bf Wild:} $X$ supports families of
arbitrarily large dimension of indecomposable non-isomorphic ACM sheaves.
\end{itemize}
\end{thm0}
{ Finite and tame cases are completely classified. Precisely we have from \cite{FP} the following list.
\vskip 0.4cm

\noindent{\bf Finite} connected $X\subset\PP^m$ are
\begin{itemize}
{\it \item{}(i) $\PP^m$
\item{}(ii)
a smooth quadric hypersurface;
\item{}(iii) a smooth rational curve
\item{}(iv) the smooth cubic scroll in $\PP^4$, isomorphic to $\PP(\O_{\PP^1}\oplus\O_{\PP^1}(1))$,  also to $\PP^2$ blown up in a point,
\item{}(v) the Veronese surface in $\PP^5$.}
\end{itemize}
{\bf Tame} connected $X\subset\PP^m$ are
\begin{itemize}
{\it \item{}(i) quadric hypersurface of corank one;
\item{}(ii) a reducible tree of rational curves;}

in cases (i), (ii) the parameter space of indecomposable non isomorphic ACM sheaves is a countable set of points,
{\it \item{}(iii)
a smooth elliptic curve or a cycle of rational curves;
\item{}(iv)
the smooth quartic scroll in $\PP^5$,  there are two kinds of them, isomorphic to $\PP^1\times\PP^1$ embedded with $\O(1,2)$, 
or isomorphic to $\PP(\O_{\PP^1}\oplus\O_{\PP^1}(2))$, see \cite[Ex. 2.19.1]{Hart}.}

in cases (iii), (iv) the parameter space of indecomposable non isomorphic ACM sheaves is a curve.
\end{itemize}

Note that the Del Pezzo surface which is the complete intersection of two quadrics in $\PP^4$ is wild, since it does not appear in the above list.
Also the Segre $3$-fold $\PP^1\times\PP^2$ (a variety of minimal degree) and $\PP^1\times\PP^1\times\PP^1$ in their minimal equivariant embedding
are wild, see \cite{CFM} for the latter.

The paper \cite{FP} by Faenzi and Pons-Llopis is the culmination of a large amount of work, beyond the two authors,  by Buchweitz, Casnati, Casanellas, Costa, Drozd, Eisenbud, Greuel, Hartshorne, Herzog, Malaspina, Mirò-Roig, Schreyer and other people that it is used in the paper and it cannot be resumed here, and we refer to the bibliography in \cite{FP} for further informations.

 Also the trichotomy has classical roots, C. Segre proved in 1886 \cite{S2} that {\it elliptic curves} are {\it tame}, a result commonly attributed to Atiyah.
For a modern exposition of this Segre result we quote  \cite[\S 5]{GO}.

\begin{rem0} 
{An updated reference} which explores the related topic of Ulrich bundles is the book \cite{CMP} by
Laura Costa , Rosa María Miró-Roig and Joan Pons-Llopis. See the nice survey \cite{Beau}. Ulrich bundles are special ACM bundles. The conjecture that every projective variety carries a Ulrich bundle is the next challenge. This is another interesting  story beyond the limits of the present paper.
\end{rem0}

\end{document}